\documentclass{article}

\usepackage{amsfonts}
\usepackage{amssymb}
\usepackage{mathrsfs}
\usepackage{amsthm}
\makeatletter
\def\@biblabel#1{#1.}
\makeatother
\newcounter{obr}[section]

\newcounter{pvv}[section]
\renewcommand{\thepvv}{\thesection.\arabic{pvv}}
\newenvironment{pv}[2][]{\begin{trivlist}\refstepcounter{pvv}%
\item[\hspace{\labelsep}\normalfont\bfseries\thepvv. #2%
  \def\tmp{#1}\ifx\tmp\empty\else{} (#1)\fi.]}%
{\end{trivlist}}
\newenvironment{df}{\begin{pv}{Definition}}{\end{pv}}
\newenvironment{theo}[1][]{\begin{pv}[#1]{Theorem}\begin{itshape}}{\end{itshape}\end{pv}}
\newenvironment{pr}{\begin{pv}{Proposition}\begin{itshape}}{\end{itshape}\end{pv}}
\newenvironment{lem}{\begin{pv}{Lemma}\begin{itshape}}{\end{itshape}\end{pv}}

\newenvironment{cor}{\begin{pv}{Corollary}\begin{itshape}}{\end{itshape}\end{pv}}

\newcommand{\rr}{\ensuremath{\mathbb{R}}}

\def\sp{\mathop{\rm sp}\nolimits}

\def\card{\mathop{\rm card}\nolimits}

\newcommand{\f}{\ensuremath{\varphi}}
 \newcommand{\ps}{\ensuremath{\psi}}
\renewcommand{\l}{\ensuremath{\lambda}}
\newcommand{\ro}{\ensuremath{\varrho}}

\newcommand{\ifff}{if, and only if,}
\newcommand{\vNa}{von Neumann algebra}

\newcommand{\Ca}{$C^\ast$-algebra}

\newcommand{\HS}{Hilbert space}

\newcommand{\st}{such that}

\newcommand{\seq}[2]{\ensuremath{#1_1, \ldots, #1_{#2}}}

\newcommand{\set}[2]{\ensuremath{\{ #1\,|\; #2\}}}

\newcommand{\Cl}[1]{\ensuremath{{\cal #1}}}

\newcommand{\un}{{\bf 1}}
\begin{document}

\begin{center}
{  \LARGE \bf Structure of associative subalgebras of Jordan operator algebras}\\[6mm]
{\large Jan Hamhalter and Ekaterina Turilova}\\[2mm]
Czech Technical University in Prague - El. Eng.\\[2mm]
Department of Mathematics, Technicka 2,\\[2mm]
 166 27 Prague 6, Czech Republic\\[2mm]
hamhalte@math.feld.cvut.cz \\[4mm]

Kazan Federal  University,\\[2mm]
 Faculty of Computational Mathematics and Cybernetics,\\[2mm]
 Kremlevskaya 18, Kazan, Russia,\\[2mm] 
e-mail: ekaterina.turilova@ksu.ru

\end{center}
\vspace{2cm}

{\scriptsize Abstract:  We show that any order isomorphism between ordered structures of associative unital JB-subalgebras of  JBW algebras is implemented naturally by a Jordan isomorphism.  Consequently, JBW algebras are determined by the structure of their associative unital JB subalgebras. Further we show that in a similar way it is possible to reconstruct Jordan structure from the order structure of associative subalgebras endowed  with an orthogonality relation. In case of   abelian subalgebras of \vNa s we show that  order isomorphisms of  the structure of abelian  $C^\ast$-subalgebras that are well behaved with respect to the structure of two by two matrices over original algebra  are implemented by  $\ast$-isomorphisms.        }

\vspace{5mm}

{\scriptsize keywords: Jordan algebras, structure of associative subalgebras}

\section{Introduction and Preliminaries}
  
 JB algebra  is a far reaching generalization of the self-adjoint part $(A, \circ)$ of a \Ca{} $(\Cl A, \cdot)$ endowed with  the product 
 \[ x\circ y=\frac 12 (xy+yx)\qquad x,y\in A\,.\]
Jordan product  is not associative in general. An  associative Jordan algebras are quite  special. More precisely, 
basic result of spectral theory says that  JB algebra $(A, \circ)$ is associative \ifff{}  it is isomorphic to the algebra $C_0(X)$ of all continuous real-valued functions on a locally compact Hausdorf space $X$ vanishing at infinity.  
 In a similar way, JBW algebras, i.e. JB algebras with preduals, are generalisation of \vNa s. 
The associative subalgebras of a given  JB algebra are mutually overlapping and, when ordered by set-theoretic inclusion, they form a partially ordered set (poset for short). This poset can be classified as a semi-lattice: any two elements admit an infimum, which is a set theoretic intersection. If two JBW algebras are isomorphic, then their corresponding posets of associative subalgebras are order isomorphic, and so we obtain an order  theoretic invariant of JBW algebras.  An interesting question arises whether, on the contrary, the poset  of associative subalgebras determine fully  structure of JBW algebras.  More precisely, in this note we shall  be mainly concentrated on  on the following problem: Let $A_1$ and $A_2$ be JBW algebras and $\Cl L(A_1)$ and $\Cl L( A_2)$ be  some structures of associative JB subalgebras of $ A_1$ and $ A_2$, respectively, ordered by set theoretic inclusion.  Let $\Cl L(A_1)\to \Cl L(A_2)$ be an order isomorphism. Is there  a unique Jordan isomorphism $\psi:A_1\to A_2$ \st{} 
\[  \f(C)=\psi(C) \qquad \mbox{ for all } C\in \Cl L(A_1)?\]     
 We say in this case that the Jordan isomorphism $\psi$ implements $\f$. The positive answer to this question means that, perhaps surprisingly, complicated functional-analytic structure of JBW algebras can be encoded into ``discrete'' algebraic ordered structure. We show that the solution  to the problem stated above depends on the type of associative subalgebras  chosen. First we show that the problem has positive answer 
if we consider associative unital JB subalgebras. (Unital subalgebra is a subalgebra containing the unit of the whole algebra.) In the sequel let $ASU(A)$ denote the poset of all unital associative JB subalgebras of a JB algebra  $A$. The following result shows that 
unital in case when $A$ is associative then the structure $ASU(A)$ determines $A$.
 This  result has been proved for abelian  subagebras in \cite{Hamhalter}. Its translation to associative Jordan algebras is straightforward.  

\begin{theo}\label{associative}
Let $A$  and  $B$  be unital associative JB algebras and $\f: ASU(A)\to ASU(B)$ be an order isomorphism. Suppose that $\dim A\not= 2$.  Then there is a unique Jordan isomorphism $\psi:  A\to B$ that  implements \f.     

\end{theo} 

In the first part of the paper we extend  this result in the following way. Suppose that $A$ is a JBW algebra without Type $I_2$ direct summand that is not isomorphic to $\rr\oplus \rr$.  Then, for each order isomorphism $\f:ASU(A)\to ASU(B),$ where $B$ is another $JBW$ algebra,  there is a unique Jordan isomorphism $\psi: A\to B$, implementing \f. In the light of counterexamples we show that the assumptions are not superfluous. 
Our result generalizes the corresponding results for abelian $C^\ast$-subalgebras of \vNa s obtained in \cite{Hamhalter} by J.Hamhalter. \\

An important types   of associative subalgebras are singly generated subalgebras. These algebras may  not have a unit in general.  It is therefore natural to consider the problem of determining whole Jordan structure by structure of  all associative subalgebras. Let us denote by $AS(A)$ the system of all associative JB subalgebras of a JBW algera $A$, again ordered by set theoretic inclusion. In this case the situation is different. Since any Jordan isomorphism preserves    the unit,  it can implement   only order isomorphism that preserves unital associative subalgebras. Therefore,  in order to recover the structure of the algebra, we have to consider some additional  structure on the set 
$AS(A).$ One of the natural possibilities seems to be  an  orthogonality relation on $AS(A)$ introduced as follows.
Let $C_1, C_2\in AS(A)$. We say that $C_1$ and $C_2$ are orthogonal
(in the symbols $C_1\perp C_2$)  if $C_1$ and $C_2$ operator commute,  and $x\circ y=0$ whenever $x\in C_1$ and $y\in C_2$. We shall show that 
more or less any transformation of 
the structure of associative subalgebras 
 that preserves  the order and orthogonality relation in both directions,  is implemented by a unique  Jordan isomorphism. This result is new also in the context of \vNa s.\\
 In case of \vNa s it is known that the order structure of abelian   subalgebras does not suffice to    recover the whole structure. The reason is a deep result   of Connes \cite{Connes},  who showed that there is a von Neumann factor $M$ that is not anti-isomorphic to itself.  It means that $M$ and its opposite algebra $M^o$  have identical structures of  all  abelian subalgebras, although they are not  $\ast$-isomorphic. (Let us recall that opposite algebra has identical linear structure and norm but the multiplication is reversed.)  In this light, the following  open problem has been formulated by A.D\"{o}ring and J.Harding \cite{DoringHarding}: What additional structure on the poset of  abelian subalgebras suffices for recovering whole \vNa{} structure? In the last part of this paper we show one possibility in this direction using the structure of two by two matrices over original algebra. \\
 
 We briefly discuss the relevance of the results obtained to foundations of quantum theory.\\

  Le us now recall a few concepts and fix the notation. (Our standard reference for Jordan algebras is \cite{Olsen}.)  JB algebra is a  real Banach algebra with the commutative product, $\circ$, that satisfies the following properties for all $a,b\in A$: (i)  $a\circ (b\circ a^2)=(a\circ b)\circ a^2$\     (ii): 
$\| a^2\| \le \|a^2+b^2\|$, $\|a\|^2=\|a^2\|$. 
JBW algera is a JB algebra that is a dual Banach space. Throughout the paper    $A$ shall denote the JBW algebra.  We write $A_1=\set{a\in A}{\|a\|\le 1}$,  $A^+=\set{a^2}{a\in A}$.  For $a\in A$, mappings $T_a, U_a : A\to A$ are defined by letting $T_a(b)=a\circ b$, $U_a(b)=2a\circ (a\circ b)-a^2 \circ b.$
It is well known that $U_a$ is positive in the sense that $U_a(A^+)\subset A^+$. By a projection we mean an idempotent of $A$. The symbol $P(A)$ shall be reserved for the set of all projections in $A$. We shall often use the fact that, for a projection $p\in A$,  $U_p(A)$ is a JBW subalgebra of $A$ with the unit $p$. For $\seq xn \in A$ we shall write $JB(\seq  xn)$ for a JB subalgebra of $A$ generated by $\seq xn$.  Two elements $a,b\in A$ are said to be operator commuting if $T_aT_b=T_bT_a$. Two subsets  $C,D\subset A$ are called operator commuting if each element of $C$ operator commutes with each element of $D$.       The center $ Z(A)$ of $A$    is the set of all elements that operator commute with each element of $A$.  If the center is trivial, then $A$ is called a factor. 
 Positive functional \f{} on $A$ is a a linear functional \st{} 
 $\f(a)\ge 0$ for all $a\in A_+$. A state is a norm one positive functional. 
We shall need the concept of the  spin factor.  Let $H_n$ be a real \HS{} of dimension 
$n\ge 2$. Let $V_n=H_n\oplus \rr$ have the norm $\|a+\l 1\|=\|a\|+\|\l\|$, $a\in H_n,$ $\l\in \rr$. Define a product in $V_n$ by 
\[ (a+\l1)\circ (b+\mu 1)= (\mu a+\l b)+(<a,b>+\l\mu)1\,,    \]
where $a,b\in H_n$, $\l,\mu\in \rr$. Then $V_n$ is a JBW algebra that is called the spin factor. It is known  that spin factors are exactly Type $I_2$ factors in structural theory of JBW algebras.\\

  Le us now recall a few concepts of theory of ordered structures.
 Suppose that $(X,\le)$ is a poset with a smallest element $\bf 0$.    
  We say that element $x\in X$ covers element $y\in X$ (in the symbol $x\triangleright y)$, if $x\ge y$, $x\not=y$,  and the following implication holds:
if $x\ge z\ge y$, then either $x=z$ or $y=z$. The elements that cover $\bf 0$ are called atoms. Let $x\in X$ and suppose that there is a finite sequence $\seq xn$ \st{} $x=x_1\triangleright x_2\triangleright \cdots \triangleright x_n\triangleright {\bf 0}$. The minimal $n$ with this property is then called the height of an element $x$.

\section{Unital associative subalgebras}

Let $\bf 1$ denote the unit of the JBW algebra $A$. In this section we shall deal with order isomorphisms of the structure $ASU(A)$ of all unital associative JB subalgebras of $A$ containing $\bf 1$,  ordered by set theoretic  inclusion. 
Let $B$ be another JBW algebra. The bijection $\f: ASU(A)\to ASU(B)$ is called an order isomorphism if it preserves the order in both directions, i.e. if
\[ C_1\subset C_2 \Longleftrightarrow \f(C_1)\subset \f(C_2) \qquad \mbox{ for all } C_1,C_2\in ASU(A)\,.\]
We say that order isomorphism \f{} is implemented by a Jordan isomorphism $\psi:A\to B$ if 
\[ \f(C)= \psi(C)\qquad \mbox{ for all } C\in ASU(A)\,.\]
The main result of this section is the following theorem:

\begin{theo}\label{theo-unital}
Let $A$ be a JBW algebra without Type $I_2$ direct summand that is not isomorphic to $\rr\oplus \rr$. Let $B$ be another JBW algebra. For any order isomorphism 
\[ \f: ASU(A)\to ASU(B)\,,\]
there is a unique Jordan isomorphism $\psi: A\to B$, implementing \f. 
\end{theo}
 
 We shall prove this theorem in series of statements.

\begin{lem}\label{dimension}
Let $X$ be a  locally compact Hausdorff space. 
If $\card X\ge k \ge 2$,  then  $C_0(X)$  contains a proper unital subalgebra of dimension at least  $k-1$.
\end{lem}

Proof: By the assumption, there are distinct points $\seq xk$ in $X$. Let us consider a subalgebra $C=\set{f\in C_0(X)}{f(x_1)=f(x_k)}$.  It is a proper unital  subalgebra. Using Uryson lemma  we argue that  $C$ contains $k-1$ functions $g_i$, ($i\le k-1$),  \st{} $g_i(x_j)=\delta_{ij}$
($i,j=1,\ldots ,k-1$). Therefore,  $\dim C\ge k-1$. \hfill $\square$\\

Let us remark that the algebra $ASU(A)$ has  a smallest element -- the space generated by \un.  Next lemma characterises the atoms  in the poset $ASU(A)$. 

\begin{lem}\label{atom}
$C\in ASU(A)$ is an atom \ifff{} there is a projection $p$ in $A$, different from {\bf 0}  and $\un$, \st{}
\[ C=\sp\{p, \un -p\}\,.\]
\end{lem}

Proof: Any two-dimensional unital associative algebra contains only one-dimensional subalgebra generated by the unit $\un$ as its proper
subalgebra and so it is an atom  in the poset $ASU(A)$.  Let us prove the reverse implication.
Suppose that $C$ is an atom  in $ASU(A)$. Then $C$ is isomorphic to $C(X)$, where $X$ is a compact Hausdorff space. If $C$ has three distinct points, then, according to
 Lemma~\ref{dimension},  $C$ contains a proper two-dimensional subalgebra. This is a contradiction. So $X=\{x,y\}$ and, in turn, $C=\sp\{p, \un -p\}$, where $p$ is a projection corresponding to the characteristic function of the set $\{x\}$.  \hfill $\square$\\

\begin{lem}\label{maximal}
Let $C$ be a maximal unital associative subalgebra of $A$ and $p\in C$ be a projection. Then 
$\dim U_p(C)=1$ \ifff{} $\dim U_p(A)=1\,.$
\end{lem}

Proof: The reverse implication  is trivial and so we show that $\dim U_p(C)=1$ implies $\dim U_p(A)=1$. Suppose, for a contradiction, that   $\dim U_p(C)=1$ and $\dim U_p(A)\ge 2$.
Then there is a positive $a$, with $0\le a \le p$,  that is not a multiple of $p$. We can write

\begin{equation}\label{star}
C= p\, \rr\oplus (\un-p)\circ C\,.
\end{equation}
As $a\in U_p(A)$, $a$  operator commutes with all elements of $U_{\un-p}(A)$ (see \cite[Lemma 2.6.3., p. 48]{Olsen}). But  $U_{\un-p}(A)$ contains $U_{\un-p}(C)=(\un-p)\circ C$,  and so, by (\ref{star}), $a$ operator commutes with $C$. In view of maximality of $C$, we obtain that $a\in C$. But then $a\in\,  \rr\,  p$, which is a contradiction.    \hfill $\square$\\

The following proposition characterises the JBW algebras that admit associative unital subalgebra that is simultaneously minimal and maximal.

\begin{pr}\label{2-maximal}
$A$ contains a two-dimensional maximal associative unital subalgebra \ifff{} $A$ is a Type $I_2$ factor or $\rr\oplus \rr$.
\end{pr}
 
 Proof: Suppose that $C$ is a two-dimensional maximal associative unital subalgebra of $A$. Then $C=\sp\{p, \un-p\}$, where $p$ is a projection. By the previous lemma 
\[ \dim U_p(A)=\dim U_{\un-p}(A)=1\,.\]
It means that $p$ and $\un-p$ are atomic projections in $A$. Let us denote by $c(p)$ and $c(\un-p)$ central supports of $p$ and $q$, respectively. As $p$ is an atom in $A$, its central cover, $c(p)$, is an atom in the center $Z(A)$.
(Let us recall that central cover $c(p)$ is the smallest central projection $z$ \st{} $z \ge p$. 
 This fact is probably known, but since we have  not been able to find appropriate quotation, we give full argument here.   
 Suppose that $z$ is a central projection \st{} $z\le c(p)$. By \cite[Lemma 4.3.5, p. 103]{Olsen},  we can compute
 \[  c(z\circ p)=z\circ  c(p)=z\,.\]
By atomicity of $p$ we conclude that either $z\circ p=0$ or $z\circ p=p$. The former case gives that $z=0$, while the latter case means that $z=c(p)$. Therefore $c(p)$ is an atom. By the same reason $c(\un-p)$ is an atomic central projection. Now we have two possibilities: 
\[ c(p)c(\un-p)=0 \qquad\mbox{ or } \qquad c(p)=c(\un -p)\,.\]
Consider first the case,   when $c(p)$ and $c(\un -p)$ are orthogonal. As
\[ \un \ge c(p)+c(\un -p) \ge p+\un -p =\un \,,\]                
 we obtain that $c(p)=p$ and $c(\un-p)=\un-p$. Whence, $A$ is isometric to $\rr\oplus \rr$.
If $c(p)=c(\un-p)$, then 
\[ c(p) \ge p +\un -p =\un\,.\]
 and so $c(p)=c(\un-p)=\un$. Atomicity of $c(p)$ in the center yields that $A$ is a factor. As $A$ contains two orthogonal projections with unit central supports that sums to $\bf 1$, $A$ has to be a factor of Type $I_2$. \\
 
 Let us now prove the reverse implication. Any Type $I_2$ factor is a spin factor $V_n=\rr\oplus H_n$, where $H_n$ is a \HS. Each nontrivial projection in $V_n$ is of the form
 \[ \frac 12(1+\mu)\,,\]
 where $\mu$ is a unit vector in $H_n$. It follows that each nontrivial projection in $V_n$ is atomic. It can be easily verified that, for a unit vector $\mu\in H_n$\,, 
 \[ \frac 12(1+\mu)\,, \frac 12(1-\mu)\]
 is a pair of orthogonal atomic projections,  whose linear span is a two-dimensional maximal unital associative subalgebra of $A$. \hfill $\square$\\
 
 Finally, we have passed to the proof of the main theorem.\\

 Proof of Theorem~\ref{theo-unital}: Suppose that $\f: ASU(A)\to ASU(B)$ is an order isomorphism. Let $X$ be a maximal associative unital subalgebra of $A$. By Proposition~\ref{2-maximal}, we have that $\dim X\ge 3$. By Theorem~\ref{associative} we see that there is a unique Jordan isomorphism $\psi_X: X\to \f(X)$ implementing  \f. We will show  the interplay between various $\psi_C$, where $C$ runs through $ASU(A)$. Let $E$ and $D$ be two different maximal elements of $ASU(A)$. We show that 
\[ \psi_E=\psi_D \qquad \mbox{ on } \ \  C=E\cap D\,.\]
If $\dim C\ge 3$, the equality above follows from uniqueness of implementing Jordan isomorphism established in Theorem~\ref{associative}. Therefore, we concentrate on the case when  $\dim C =2$, which means that 
\[ C=\sp\{p, \un-p\}\,, \ \  \mbox{where } p\in P(A)\,.\]

It can be verified easily that if $\dim U_p(A)=\dim U_{\un-p}(A)=1$, then  $C$ is a maximal associative subalgebra. 
Indeed, if $x\in A$ operator commutes with $C$,  then $p\circ x=U_p(x)=\l_1 p$ for some $\l_1\in \rr$. Similarly,  $(\un-p)\circ x =U_{\un-p}(x)=\l_2 p$ for some $\l_2\in \rr$.  It implies that $x=p\circ x+(\un-p)\circ x =\l_1 p+\l_2 p\in C.$ \\
However, it would mean, by Proposition~\ref{2-maximal},  that $A$ is either Type $I_2$  or $\rr\oplus \rr$, which is excluded by the assumption.  So it cannot happen that $p$ and $\un-p$ are simultaneously one-dimensional projections. Consider first the case when  
\[ \dim U_p(A)=1 \qquad  \mbox{ and  } \ \  \dim U_{1-p}(A)\ge 2\,.\]
By Lemma~\ref{atom} 
\[ \f(C)=\sp\{q, \un-q\}\,\]
for some $q\in P(B)$. Obviously,
\[ \psi_D(C)=\psi_E(C)=\sp\{q, \un -q\}\,.\]
As $\dim U_p(D)=1$, we infer that 
 \[ \dim U_{\psi_D(p)}(\psi_D(D))=1\,.\]                     
In virtue of Lemma~\ref{maximal}, we have that 
\[\dim U_{\psi_D(p)}(B)=1\,.\]
We shall prove further that 
 \[ \dim U_{\un-\psi_D(p)}(B)\ge 2\,.\]
 Assume the opposite and try to reach a contradiction.
If $\dim U_{\un-\psi_D(p)}(B)=1\,,$ then
$\psi_D(C)\,  (=\f(C))$  is a maximal associative subalgebra of $B$ (see the arguments before). As 
\f{} is isomorphism, the same must be true of $C$. It would mean that $C=D=E$, contradicting $E\not= D$. Suppose, without loss of generality, that $q=\psi_D(p)$. By applying the same arguments to $\psi_E$ we have that $\psi_E(p)=q,$ and so $\psi_D=\psi_E$ on $C$.\\

  As the second case, suppose that $\dim U_p(A)\ge 2$ and $\dim U_{1-p}(A)\ge 2$. Employing
  Lemma~\ref{maximal}, there exist elements $x\in D, y\in E$ \st{} $x\le p, y\le \un -p$ and \st{} the sets $\{p, x\}$ and $\{\un-p, y\}$ are linearly independent. Set 
\[ V=JB(\un,p,x)\,, \quad W=JB(\un,p, y)\,,\quad  U=JB(\un, p, x,y)\,.\]
These algebras are associative. It holds that $\dim V = \dim W\ge 3$. By Theorem~\ref{associative}  we have
\begin{eqnarray}
\psi_V(p) &=& \psi_D(p)\\
\psi_W(p) &=& \psi_E(p)\,.
\end{eqnarray}  
But the inclusions $V,W\subset U$ entail
\[ \psi_V(p)=\psi_W(p)=\psi_U(p)\,.\]
We can conclude again that 
  \[ \psi_D = \psi_E \qquad \mbox{ on }  C\,.\]
  Let us now define the map $\psi$. If $x\in A$, then $x$ is contained in some maximal associative unital JB subalgebra $X$. Put
\[ \psi(x)=\psi_X(x)\,.  \]  
By the previous reasoning this definition does not depend on the choice of $X$. As $\psi_X$ is a Jordan isomorphism, we see that $\psi$ is an isometry. The restriction of $\psi$ to the projection lattice $P(A)$ gives a finitely additive measure, which means that
\[\psi(p+q)=\psi(p)+\psi(q)\,,\]
whenever $p$ and $q$ are orthogonal projections. Let us consider a state $f\in B^\ast$. Then  
 $f\circ \psi$ is a nonnegative scalar measure. By deep Gleason type theorem for JBW algebras proved by Bunce and Wright in \cite{BW1,BW2}, $f\circ \psi$ extends uniquely to a state $\psi_f$ on $A$. (For this argument the assumption on Type $I_2$ factor is necessary.)  Take $x,y\in A$. Then
 \[ f(\psi(x+y))=\psi_f(x+y)=\psi_f(x)+\psi_f(y)=f(\psi(x)+\psi(y))\,.\]
 As the set of states is separating, we see that
 \[ \psi(x+y)=\psi(x)+\psi(y)\,,\] 
 establishing the linearity of $\psi$. It follows immediately from the definition of $\psi$,  
 that $\psi$ is a Jordan homomorphism that implements \f. Applying the same argument to $\f^{-1}$, we see that $\psi$ is a Jordan surjective isomorphism. \\
 
    Finally, let us establish the uniqueness of $\psi$. Suppose that \f{} is implemented by two Jordan isomorphisms $\psi_1, \psi_2: A\to B$. Then $\ro=\psi_1\circ \psi_2^{-1}$ is a Jordan isomorphism that implements the identity on $ASU(A)$. For any $x\in A$ take a maximal
associative unital subalgebra $C$ of $A$ containing $x$. By assumption, $\dim C\ge 3$, and so 
$ \ro(x)=x$, for all $x\in C$ by Theorem~\ref{associative}. In other words, $\ro$ is an identity and so $\psi_1=\psi_2$. \hfill $\square$\\

     Let us remark that the assumption on the absence of Type $I_2$ factor is not superfluous  for establishing the uniqueness of implementing Jordan isomorphism. Indeed,
consider spin factor $V_n=\rr\oplus H_n$. Then, it can be verified easily that all associative unital subalgebras (except for whole algebra and the algebra generated by the unit) are of the form
\[ A_\xi= \sp\biggl\{\frac 12 (1+\xi), \frac 12(1-\xi)\biggr\}\,,  \]                                  where $\xi$ is a unit vector.  System $(A_\xi)_\xi$, where $\xi$ runs through the unit ball of $H_n$ is a family of  subalgebras covering whole $V$ and  intersecting each other in the subalgebra $\sp\{\un\}$. Let $\psi_\xi$ be a Jordan automorphism of $A_\xi$. The union of 
the graphs of $\psi_\xi's$ gives a Jordan  automorphism of $A$. By this way we can get many Jordan automorphisms of $V$ implementing the identical order isomorphism of $ASU(V_n)$.

   \section{Non unital associative subalgebras}

   This section is devoted to reconstruction of JBW algebra from its structure of all nonunital associative subalgebras.  Let us denote by  $AS(A)$ the structure of all associative JB subalgebras of $A$.    As we remarked already in the introduction, the order isomorphism  between posets of associative subalgebras may not be implemented by any Jordan isomorphism. For example, let us take an algebra $\rr\oplus \rr$. Then   $AS(A)$ consists of  the whole algebra (the largest element), the zero algebra (the smallest element), and then of three atoms. Any permutation of the atoms  that leave other algebras unchanged  is  an order isomorphism. One of the atoms is the algebra generated by the unit. This atom is invariant with respect to any  isomorphism implemented  by 
a Jordan isomorphism. So there are order automorphisms of $AS(A)$ that are not implemented by Jordan automorphisms. In the light of it, we shall introduce another structure on $AS(A)$ that will complete the information contained in the order.    We say that two algebras $C$ and $D$  in $AS(A)$ are orthogonal (in symbols $C\perp D$) if    $x\circ y=0$   for all $x\in C$ and $y\in D$. Let us note that orthogonality  implies that $C$ and $D$ operator commute.  Indeed, let us take positive elements $x\le {\bf 1}$ and $y\le {\bf 1}$ from $C$ and $D$, respectively. Since for any polynomial $p$ with $p(0)=0$ we have that $p(x)\in C$, we infer  that $p(x)\circ y=0$. Consequently, using function calculus we can  see that
for any continuous function $f$, with $f(0)=0$, we have $f(x)\circ y=0$. As the range projection $r(x)$ of $x$  is a strong operator limit of elements of the form $f(x)$ above,
we conclude that $r(x)\circ y=0$. By symmetry we see that $r(x)$ and $r(y)$ are orthogonal projections. By the well known result $U_{r(x)}(A)$ and $U_{r(y)}(A)$ operator commute.
As $x\in U_{r(x)}(A)$ and $y\in U_{r(y)}(A)$, we obtain that $x$ and $y$ operator commute. Since any element in a JB algebra is a difference of two positive elements, we have shown that $C$ and $D$ operator commute.

Next theorem is the main result of this section.

\begin{theo}\label{orthogonality}
Let $A$ and $B$ be JBW-algebras \st{} $A$  has no Type $I_2$ direct summand. Let $\f: (AS(A), \subset, \perp)\to (AS(B), \subset, \perp)$ 
be an order isomorphism preserving the orthogonality, $\perp$,  in both directions. Then there is a unique Jordan isomorphism  
\[ \psi: A \to B\]
\st{}
 \[ \f(C)=\psi(C) \qquad \mbox{ for all } C\in AS(A)\,.\] 
 \end{theo}

Proof  will be divided into a few next steps. 

\begin{pr}\label{height}
Let $A$ be a JBW algebra. 
The following conditions are equivalent:
\begin{enumerate}
\item $X\in AS(A)$ has dimension $k$
\item $X=\sp\{ \seq pk\}$, where $\seq pk$ are orthogonal nonzero projections. 
\item There is a chain in $AS(A)$ \st{}
\[ X \triangleright  X_  {k-1}\triangleright \cdots \triangleright X_1\triangleright \{0\}\,.\]                               
\end{enumerate}

\end{pr}
 
 Proof: $(i)\Rightarrow (ii)$ Each associative JB algebra is isomorphic to the algebra $C_0(X)$ of real valued continuous functions on a locally compact  Hausdorff topological space $K$ that vanish at infinity. If $\dim C_0(X)=k$, 
 then $X$ must contain precisely $k$ points, i.e. $K=\{x_1, \ldots, x_k\}$.  Characteristic functions, $p_i$,  of singletons $\{x_i\}$ give the desired system of projections.  \\
 
 $(ii) \Rightarrow (iii)$ It is immediate by setting
 \[ X_{k-1}=\sp\{p_2, \cdots, p_k\}, \ldots, X_1=\sp\{p_1\}\,.\]
 
 $(iii)\Rightarrow (ii)$ We shall proceed by induction on $k$. Suppose that $X=C_0(K)$ is an atom in $AS(A)$. We shall prove that $K$ is singleton. Indeed, assume  $K$ has two distinct points, say $x$ and  $y$. They  can be separated by a continuous function $f$ in $C_0(K)$ in the sense that $f(x)=1, f(y)=0$. The algebra $JB(f)$ is an associative nonzero subalgebra consisting of  functions vanishing at $y$. It means that this algebra is a proper nonzero subalgebra of $X$, which contradicts  the fact that $X$ is an atom. 

  Suppose now that $X\triangleright X_{k-1}$, where $\dim X_k=k-1$. We shall prove that $\dim X_k=k$. Suppose, for a contradiction, that $\dim X\ge k+1$,  which forces that $\card K\ge k+1$. Let us represent $X$ as a space $C_0(K)$ again. Then $X_{k-1}$ is by induction hypothesis a linear span of projections $\seq p{k-1}$, each corresponding to   a clopen subset $\seq Z{k-1}$ of $K$.  Then $\seq Z{k-1}, Z_k=K\setminus\bigcup_{i=1}^{k-1} Z_i$ is a covering of $K$ by disjoint clopen sets. As $\card K \ge k+1$, at least one of these sets, say $Z_1$, contains two distinct points, say $x$ and $y$.  Using basic properties of  locally compact spaces, there 
are function $f,g$ in $C_0(X)$ \st{} $f(x)= 1$, $f(y)=0$, $g(x)=0$, $g(y)=1$,  and \st{} $f$ and $g$ vanish outside $Z_1$. It is then clear that the algebra $JB(f, g, p_2, \ldots, p_k)$, where $p_k$ is a projection corresponding to $Z_k$,  is  a proper subalgebra of  $X$ that contains 
$X_{k-1}$ as its proper subalgebra. This is a contradiction. 
\hfill $\square$

\begin{cor}\label{atoms}
There is a one-to-one correspondence between atoms in $AS(A)$, and nonzero projections in $A$ \st{} each atom $C\in AS(A)$ is of the form $C=\sp\{p\}$, where $p\in P(A)\setminus\{0\}$.    
\end{cor}

\begin{lem}\label{decomposition}
Let $p$ and $q$ be  nonzero orthogonal projections and  $z$ be a nonzero projection different from $p$ and $q$.  Then
\[ z=p+q\,,\]
\ifff{} the set  $\{p,q, z\}$ is a 2-dimensional asociative subalgebra.   
\end{lem}

Proof: If $ \{p, q, z\}$ is a two-dimensional associative algebra, then  $z$ must be equal to $p+q$. Indeed, if $z\not= p+q$, then  for one of the remaining  projections, say $p$, we have $z\perp p$. But then $q=z+p$,  which is in contradiction with $q\perp p$.   The reverse implication is obvious. \hfill $\square$\\

Proof of the Theorem~\ref{orthogonality}: Suppose that $\f:AS(A)\to AS(B)$ is an order isomorphism preserving the orthogonality in both directions. Since \f{} preserves faithfully one-dimensional subalgebras (i.e. atoms), we can define a map
\[ \psi: P(A)\setminus\{0\}\to  P(B)\setminus\{0\}\]
 by equation 
 \[ \sp\{\psi(p)\}=\f(\sp\{p\})\,.\]    
 Then, according to Corollary~\ref{atoms}, $\psi$ is a bijective map of $P(A)$ onto $P(B)$. As two nonzero projections $p,q$ are orthogonal \ifff{} $\sp\{p\}$ and $\sp\{q\}$ are orthogonal algebras, 
we see that $\psi$ preserves orthogonality of projections.  Suppose now that
\[ z=p+q\,,  \] 
where $z,p,q$ are nonzero projections. Thanks to the fact that \f{} preserves two-dimensional subalgebras (Propositon~\ref{height}), we have 
\[ \f(\sp\{p, q\})=\sp\{\psi(p), \psi(q)\}\,,                  \]
 and by Lemma~\ref{decomposition} we see that 
 \[ \psi(z)=\psi(p+q)=\psi(p)+\psi(q)\,.\]
 In other words, defining $\ps(0)=0$,  $\psi: P(A)\to P(B)$ is a bounded finitely additive vector measure. 
Using now Gleason type theorem and the same arguments   as in the proof of Theorem~\ref{theo-unital},   we  infer that  $\psi$ extends to a unique bounded linear map
between $A$ and $B$. We shall denote this map by the same symbol. As $\psi$ preserves the projections, it is Jordan homomorphism. Applying the same argument to the inverse, $\f^{-1}$, we deduce that $\psi$ is 
a Jordan isomorphism. Let us now take  an arbitrary JB subalgebra $C\in AS(A)$. Fix an element $0\le x \le 1$ in $C$. Then 
\[ x = \sum_n \frac 1{2^n} p_n\,,\]
where $(p_n)$ is a sequence of projections in $C$. As $\psi(p_n)\in \f(C)$ for all $n$, we conclude that 
\[   \psi(x)=\sum_n \frac 1{2^n}\psi(p_n)\in \f(C)\,.  \]  
Hence, $\psi(C)\subset \f(C)$. Applying the same argument to the inverse map $\f^{-1}$ we see that 
\[ \psi(C)=\f(C)\,.\]

Uniqueness of $\psi$ is a consequence of the fact that Jordan map is
uniquely determined by its values on the projections. 
\hfill $\square$\\

Let us notice that, by the arguments in the proof, it is, in fact, possible to recover Jordan structure from the structure of elements of height two in $AS(A)$.

\section{Recovering full $C^\ast$-structure}

In this section we shall make an initial step to finding structure on abelian subalgebras  that contain  all information on whole 
algebra already  and not only on the Jordan part of the associative product. 
One of the natural possibilities  is to consider the structure of two by two matrices over original algebra. As it is well known, for example,  Jordan maps  that are well behaved with respect to this structure, are already $\ast$-homomorphisms. This motivates the following definition that can be viewed as order theoretic  version of the amplification for maps.   First,  some notation will be useful.   In the sequel let $\Cl A$ be a \vNa{} with self-adjoint part $ A$. $M_2(\Cl A)$ will denote the algebra of two by two  matrices over $\Cl A$.   We will identify it  with the tensor product $M_2\otimes \Cl A$, where $M_2$ is the algebra of  two by two complex matrices. The algebra $\Cl A$ will be embedded into  $M_2(\Cl A)$  diagonally as $\un \otimes \Cl A$. By $ Abel(\Cl A)$   we shall mean the structure of abelian  subalgebras  of $\Cl A$ ordered by set inclusion. The orthogonality relation, $\perp$, on $Abel(A)$ will be naturally extended from the corresponding structure $AS(A)$: We say
that $C,D\in Abel(\Cl A)$ are orthogonal (in symbols $C\perp D$)  if $xy=0$ for all $x\in C$ and $y\in D$. 

\begin{df}
Let $\Cl A$ and $\Cl B$ be \vNa s.  
Let $\f: Abel(A) \to Abel(B) $
be an order isomorphism preserving the orthogonality in both directions. We say that $\omega_: Abel(M_2(\Cl A))\to Abel(M_2(\Cl B))$ is an {\em amplification} of \f{} if the following conditions are satisfied. 
\begin{enumerate}
\item $\omega$ is an order isomorphism preserving orthogonality in both directions. 
\item For all abelian subalgebras $C\in Abel(\Cl A)$ and $D\in Abel(\Cl  B)$ we have
\[ \omega(\Cl C\otimes \Cl D)= \Cl C\otimes \f(\Cl D)\,.\]

\end{enumerate}

\end{df}

\begin{pr}
Suppose that $\Cl A$ has no direct summand of Type $I_2$ and $\Cl B$ is another \vNa.  Let $\f: Abel(\Cl A)\to Abel(\Cl B)$ be  an order isomorphism preserving the orthogonality in both directions that has an amplification $\omega$. Then \f{} is implemented by a $\ast$-isomorphism between $\Cl A$ and $\Cl B$. 
\end{pr}

Proof: By our results from the previous section (rewritten easily into the language of \vNa s),  $\omega$ is induced by a Jordan $\ast$-isomorphism $\Psi:M_2(\Cl A)\to M_2(\Cl B)$,  and \f{} is induced by a Jordan $\ast$-isomorphism $\psi: \Cl A\to \Cl B$. \\

Let $p$ and $q$ be projections in $\Cl A$ and $\Cl B$, respectively. Then
\[ \sp\{p\}\otimes \sp\{q\}=\sp\{p\otimes q\}                 \]  
is an abelian subalgebra which is mapped by $\omega$ onto abelian subalgebra 
\[  \sp\{p\}\otimes \sp\{q'\}\,, \qquad \mbox{ where } q'\in P(\Cl B)\,.\]
Therefore $\Psi(p\otimes q)\in \sp\{p\otimes q'\}$ and taking into   account that $\Psi$ is a Jordan map (and so preserving projections), we infer that 
\begin{equation}\label{alfa}
\Psi(p\otimes q)=p\otimes q'\,.
\end{equation}
By the spectral theorem, we obtain from    (\ref{alfa})  that for all $x\in \Cl A$, $y\in \Cl B$ we have 
\[ \Psi(x\otimes y)= x\otimes z\,,  \]
 where $z$ is some element of $\Cl A$. In other words,
 \[ \Psi= id \otimes \ro\,,\]
 where $\ro: \Cl A \to \Cl B$ is a linear map. Because \ro{} preserves projections, it is a Jordan map.  As $\un\otimes \Cl A$ is identified with $\Cl A$, we have that $\ro=\psi$ by uniqueness of a Jordan map implementing \f. So $\psi$ is a Jordan map whose amplification is again Jordan. Consequently,  $\psi${} is 2-positive. By the result of Choi \cite[Cor. 3.2]{Choi},  any two-positive unital Jordan map is a $\ast$-homomorphism. This concludes the proof.  \hfill $\square$\\

Let us remark that in mathematical foundations of quantum theory,  the system of observables is given by a JB algebra or a JBW algebra. This model extends the classical system in which observables are given by elements of an associative JB algebra, i.e. by  the  algebra of continuous functions on some phase space. From this point of view  
the quantum system is a union of mutually overlapping and interacting classical systems.  It seems to be natural to capture the structure  of classical subsystems by the structure of associative 
subalgebras considered in this paper. In this context, our results may be viewed  as possibility to reconstruct global quantum system and its symmetries (given by Jordan isomorphisms) from classical subsystems. 
That is why, apart from its own mathematical interest, our study has been motivated by extensive research in mathematical foundations of
quantum theory (including topos theory) \cite{DoringIsham, Landsman}.

 \end{document}

  \end{document}